\documentclass[11pt]{article}
\RequirePackage[amsthm,amsmath]{imsart}

\startlocaldefs
\usepackage{amsmath}
\usepackage{amssymb}
\usepackage{amsfonts}
\usepackage{graphicx}

\usepackage{amsthm}
\newtheorem{theo}{Theorem}

\newtheorem{col}{Corollary}
\newtheorem{lem}{Lemma}

\endlocaldefs

\begin{document}

\begin{frontmatter}

\title{A Tight Bound of Tail Probabilities for a Discrete-time Martingale with Uniformly Bounded Jumps}
\runtitle{A tight bound for a discrete-time martingale}

\author{\fnms{Go} \snm{Kato}\thanksref{t1}\ead[label=e1]{go.kato.gm@hco.ntt.co.jp}}
\thankstext{t1}{This work is supported by JSPS KAKENHI Grant Number JP17K05591}

\runauthor{G. Kato.}

\affiliation{NTT Communication Science Laboratories, NTT Corporation}
\address{
3-1,Morinosato Wakamiya Atsugi-Shi, Kanagawa, 243-0198, Japan\\
\printead{e1}}

\begin{abstract}
We investigate the properties of a discrete-time martingale $\{X_m\}_{m\in \mathbb Z_{\geq 0}}$, where all differences between adjacent random variables are limited to be not more than a constant as a promise. In this situation, it is known that the Azuma-Hoeffding inequality holds, which gives an upper bound of a probability for exceptional events. The inequality gives a simple form of the upper bound, and it has been utilized for many investigations. However, the inequality is not tight. We give an explicit expression of a tight upper bound, and we show that it and the bound obtained from the Azuma-Hoeffding inequality have different asymptotic behaviors.
\end{abstract}



\end{frontmatter}

\section{Introduction and preliminary}
Providing concentration inequalities is one of the ultimate goals of probability theory. Actually, many types of such inequalities~\cite{B24,C52,B62,A67,T67,C81,H63} have been derived and applied in a huge number of fields where probabilistic events occur. However, almost all of them do not give tight bounds. Therefore, much effort has been spent in obtaining tighter inequalities~\cite{M00,S00,BL03,KR05,KR08} to provide benefits for those fields, but it is still difficult to get a  tight bound itself. The Azuma-Hoeffding inequality~\cite{A67} is one of the famous and widely used concentration inequalities. It gives an upper bound of a tail probability for a discrete-time martingale with bounded jumps, but the bound is not tight, unfortunately. We explicitly present a tight bound of the tail probability in the case where all the jumps are bounded by a single constant. The derived expression enables us to compare the tight bound with that given from the Azuma-Hoeffding inequality.

We first show the Azuma-Hoeffding inequality strictly to make this manuscript self-contained, though it is well known. We introduce a probabilistic space $(\Omega,\mathcal F,\mathbb P)$ and a martingale $\{X_n'\}_{n \in \mathbb Z_{\geq 0}}$, where all differences between adjacent random variables are limited. Precisely speaking, we consider a series of random variables which satisfies the relations
\begin{eqnarray}
\mathbb P \left(\left|X_{m+1}'-X_{m}'\right|\leq c_{m+1}\right)&=&1,
\label{eq:main_assumption_0_1}
\\
\mathbb E\left(X_{m+1}'\big|\mathcal F'_m\right)&=& X_m' \quad\makebox{\rm almost surely (a.s.)},
\label{eq:main_assumption_0_2}
\\
\mathbb E\left(X_{m}'\big|\mathcal F'_m\right)&=& X_m'
\label{eq:main_assumption_0_3}
\end{eqnarray}
for any $m\in \mathbb Z_{\geq 0}$, where $\{c_n\}_{n\in \mathbb Z_{\geq 0}}$ is an arbitrary series of real numbers, and $\{\mathcal F'_n\}_{n\in \mathbb Z_{\geq 0}}$ is a filtration, i.e., $\mathcal F'_n$ is a $\sigma$-algebra which satisfies $\mathcal F'_{n}\subset \mathcal F'_{n+1}\subset \mathcal F$ for any $n\in\mathbb Z_{\geq 0}$. The first relation, (\ref{eq:main_assumption_0_1}), is the bounded difference condition, and the other relations, (\ref{eq:main_assumption_0_2}), (\ref{eq:main_assumption_0_3}), are the conditions for $\{X_n'\}_{n\in\mathbb Z_{\geq 0}}$ to be a martingale. In this case, the relations
\begin{eqnarray}
\mathbb P\left( X_m'-X_0' \geq  x\right)&\leq & \exp({-\frac {x^2}{2\sum_{n=1}^mc_n^2}}),\\
\mathbb P\left( \left|X_m'-X_0'\right| \geq  x\right)&\leq & 2\exp({-\frac {x^2}{2\sum_{n=1}^mc_n^2}})
\end{eqnarray}
hold for non-negative $x\in\mathbb R$ and $m\in\mathbb Z_{\geq 1}$. This is the Azuma-Hoeffding inequality.

In this paper, we consider the special case where the variables $c_n$ take the same value $c$ for any $n\in\mathbb Z_{\geq 0}$, and present tight upper bounds of $\mathbb P\left( X_m'-X_0' \geq xc\right)$ and $\mathbb P\left( \left|X_m'-X_0'\right| \geq xc\right)$ in the case where $x$ is an integer. The rest of this paper is organized as follows. In section \ref{sec:main}, we explicitly explain our main result. In sections \ref{sec:lem1} and \ref{sec:lem2}, we prove lemmas which directly give our main result. In section \ref{sec:comp}, our tight bound is numerically compared with the bound given from the Azuma-Hoeffding inequality. The last section is devoted to a conclusion. 

\section{Main claim}
\label{sec:main}
Our main claim in this paper is as follows:
\begin{theo}
We assume $x,y\in\mathbb Z_{\geq 0}$ and that random variables $\{X_n\}_{n\in \mathbb Z_{\geq 0}}$ satisfy the relations
\begin{eqnarray}
\mathbb P\left(\left|X_{m+1}-X_{m}\right|\leq c\right)&=&1,
\label{eq:main_assumption_1}
\\
\mathbb E\left(X_{m+1}|\mathcal F_{m}\right)&=& X_{m}\quad\makebox{\rm a.s.},
\label{eq:main_assumption_2}
\\
\mathbb E\left(X_{m}|\mathcal F_{m}\right)&=& X_{m}
\label{eq:main_assumption_3}
\end{eqnarray}
for any $m\in\mathbb Z_{\geq 0}$, where $c$ and $\{\mathcal F_n\}_{n\in\mathcal Z_{\geq 0}}$ are an appropriate positive constant number and a filtration, respectively. In this case, the inequality
\begin{eqnarray}
\mathbb P\left( X_m-X_0 \geq xc \lor  -yc \;\geq\; X_m-X_0\right)&\leq & G(x,y,m)
\label{eq:main_value}
\end{eqnarray}
holds, where $G(x,y,m)$ for $x,y,m\in \mathbb Z_{\geq 0}$ is defined below. Furthermore, the bound is tight. Precisely speaking, for given $x,y\in\mathbb Z_{\geq 0}$ case, there are random variables $X_0$, $X_1$,$\cdots$ for which the left value in the relation (\ref{eq:main_value}) is equal to $G(x,y,m)$ for any $m\in\mathbb Z_{\geq 0}$.

The definition of $G(x,y,m)$ is as follows: When $x=0$ and $y=0$, $G(x,y,m)$ is defined to be $1$. In the other cases,
\begin{eqnarray}
&&
G(x,y,m)
\nonumber\\
&:=&
\sum_{n=0}^{\lfloor\frac{m}{2(x+y)}\rfloor}
2 I_b\bigl(\lfloor \frac{m-x}2\rfloor -(x+y)n,2\lfloor \frac{m-x}2\rfloor +x+2\bigr)
\nonumber\\&&\quad\quad{}
-2 I_b\bigl(\lfloor \frac{m-x}2\rfloor-y -(x+y)n,2\lfloor \frac{m-x}2\rfloor +x+2\bigr)
\nonumber\\&&\quad\quad{}
+2 I_b\bigl(\lfloor \frac{m-y}2\rfloor -(x+y)n,2\lfloor \frac{m-y}2\rfloor +y+2\bigr)
\nonumber\\&&\quad\quad{}
-2 I_b\bigl(\lfloor \frac{m-y}2\rfloor-x -(x+y)n,2\lfloor \frac{m-y}2\rfloor +y+2\bigr)
\label{eq:main_func}
\end{eqnarray}
where $I_b$ is a cumulative distribution function of a binomial-like distribution
\begin{eqnarray}
I_b\bigl(n,m\bigr)&:=&\sum_{z=0}^n\frac{1+(-1)^{n-z}}{2^m}\frac{m!}{z!(m-z)!}.
\label{eq:qumlant}
\end{eqnarray}
\end{theo}
Here and hereafter, we use the following notations: $\lfloor a\rfloor$ for $a\in\mathbb R$ is the floor function, i.e. the largest integer which is not larger than $a$, and $\sum_{z=x}^yf(z)$ for $y<x$ is considered to be $0$, which indicates that $I_b\bigl(-n,m\bigr)$ for $n,m\in \mathbb Z_{\geq1}$ is equal to $0$, for example.

The upper bound $G(x,y,m)$ has a closed form but a somewhat complicated one. Therefore, we give a simple bound as a corollary.

\begin{col}
For any non-negative integers $x,y,m\in\mathbb Z_{\geq 0}$, any real positive number $c\in \mathbb R$ and any martingale $\{X_n\}_{n\in\mathbb Z_{\geq 0}}$ which satisfies condition (\ref{eq:main_assumption_1}), the relation
\begin{eqnarray}
&&\mathbb P\left( X_m-X_0 \geq xc \lor  -yc \;\geq\; X_m-X_0\right)
\nonumber\\
&\leq & 
2 I_b\bigl(\lfloor \frac{m-x}2\rfloor ,2\lfloor \frac{m-x}2\rfloor +x+2\bigr)
+2 I_b\bigl(\lfloor \frac{m-y}2\rfloor ,2\lfloor \frac{m-y}2\rfloor +y+2\bigr)
\nonumber\\
\end{eqnarray}
holds.
\end{col}
This corollary can be obtained from the trivial relation $\mathbb P\left(A\lor B\right)\leq \mathbb P\left(A\right)+\mathbb P\left(B\right)$ and the explicit expression $G(x,+\infty,m)=G(+\infty,x,m)=2 I_b\bigl(\lfloor \frac{m-x}2\rfloor ,2\lfloor \frac{m-x}2\rfloor +x+2\bigr)$, which is directly derived from the definition of $G$. This is not a tight bound in a part of the region $x,y,m\in\mathbb Z_{\geq 0}$. However, it is not such a slack one (see section \ref{sec:comp})

To prove theorem 1, it is enough to prove it in the case of $c=1$ since relations (\ref{eq:main_assumption_1})$\sim$(\ref{eq:main_value}) become those for $c=1$ by replacing $X_n$ with $c X_n$ for $n\in\mathbb Z_{\geq 0}$. Therefore, in the followings, we treat only the case of $c=1$ without loss of generality. We divide the claim of the theorem into two parts: 1) The inequality (\ref{eq:main_value}) holds. 2) The bound given by the inequality (\ref{eq:main_value}) is tight. The following two lemmas are sufficient conditions for each claim:

\begin{lem}
\label{lem:uppper_bound}
We consider a series of random variables $\{X_n\}_{n\in \mathbb Z_{\geq 0}}$, which satisfies assumptions (\ref{eq:main_assumption_1}), (\ref{eq:main_assumption_2}) and (\ref{eq:main_assumption_3}) for $c=1$, and random variables $X,Y$ such that
\begin{eqnarray}
\mathbb E\left(X\big|\mathcal F_0\right)&=&X,
\label{lemma:upper_ass_1}
\\
\mathbb E\left(Y\big|\mathcal F_0\right)&=&Y,
\label{lemma:upper_ass_2}
\\
X+Y&\in& \mathbb Z_{\geq 0}.
\end{eqnarray}
The relation
\begin{eqnarray}
\mathbb P\left( X_m- X_0 \geq X\lor  -Y\;\geq\;  X_m- X_0\big|\mathcal F_0\right)
 &\leq &H_{X+Y,m}(X-Y),
\label{lemma:upper}
\end{eqnarray}
holds for any $m\in\mathbb Z_{\geq 0}$ {\rm (a.s.)}. Here, $H_{n,m}(t)$ is a piecewise linear and connective function that connects the discrete function $G(\frac12(n+t),\frac12(n-t),m)$ with respect to $t\in\mathbb R$ for any $n,m\in\mathbb Z_{\geq 0}$ in the region $|t|\leq n$, and is defined as 
\begin{eqnarray}
H_{n,m}(t)&:=&1
\label{eq:hatG(x<0,y,m)}
\end{eqnarray}
in the other regions.
\end{lem} 
This is a sufficient condition of the inequality (\ref{eq:main_value}) since both sides of Eq.(\ref{lemma:upper}) become those of Eq. (\ref{eq:main_value}) for $c=1$ by taking the average of these values in the case where both $X$ and $Y$ are constant non-negative integers $x$ and $y$, respectively.

Note that, the function $H_{n,m}(t)$ can be written explicitly as
\begin{eqnarray}
H_{n,m}(t)&=&
(1-\frac{n+t}2+z)G(z ,n-z,m)
\nonumber\\&&{}
+(\frac{n+t}2-z)G(z+1 ,n-z-1,m)
\nonumber\\
z&:=&\lfloor \frac{n+t}2\rfloor
\label{eq:hatG(x>0,y>0,m)}
\end{eqnarray}
for any $n,m\in\mathbb Z_{\geq 0}$ and $|t|< n$. The connectivity of $H_{n,m}(t)$ at the point $|t|=n$ can be checked from the facts that $G(n,0,m)=G(0,n,m)=H_{n,m}(n-0)=H_{n,m}(-n+0)=1$ (see Appendix \ref{sec:G(x,y=0,n)} ) and $H_{n,m}(n+0)=H_{n,m}(-n-0)=1$ from definition (\ref{eq:hatG(x<0,y,m)})

The tightness is given from the following lemma constructively: 

\begin{lem}\label{lem:probability_distribution} 
For given non-negative integers $x,y\in \mathbb Z_{\geq 0}$, we construct random variables $Y^{(x,y)}_1,Y^{(x,y)}_2,\cdots$, which are successively and probabilistically decided as follows: The candidates of these random variables are only $1$, $0$, and $-1$. If $\sum_{n=1}^{m-1}Y^{(x,y)}_{n}$ is in the region $\{z|x>z>-y\}$, $Y^{(x,y)}_m$ is equal to $\pm 1$ with probability $\frac 12$. In other cases, $Y^{(x,y)}_m$ is equal to $0$ with probability $1$. In the case of $X_m:=\sum_{n=1}^{m}Y^{(x,y)}_{n}$ for $m\in \mathbb Z_{\geq 0}$, assumptions (\ref{eq:main_assumption_1}), (\ref{eq:main_assumption_2}) and (\ref{eq:main_assumption_3}) for $c=1$ are satisfied with an appropriate filtration, and the value $\mathbb P\left( X_m-X_0 \geq x\lor -y\;\geq\; X_m-X_0\right)$ is equal to $G(x,y,m)$ for any $m\in Z_{\geq 0}$.
\end{lem} 

In the next two sections, we prove the lemmas.

\section{Proof of lemma 1}
\label{sec:lem1}
We prove lemma \ref{lem:uppper_bound} by the mathematical induction for $m$.

When $m$=0, we can check relation (\ref{lemma:upper}) easily from 
\begin{eqnarray}
\mathbb P\left(0\geq x\;\lor\;-y\geq 0\right)=\left\{
\begin{array}{ll}
1&\makebox{when $x\leq 0$ or $y\leq 0$},\\
0&\makebox{when $x> 0$ and  $y> 0$},
\end{array}
\right.
\\
H_{n,0}(t)\left\{
\begin{array}{ll}
=1&\makebox{when $|t|\geq n$},
\\
\geq 0&\makebox{when $t<n$}.
\end{array}
\right.
\end{eqnarray}
The first relation is trivial. The first line in the second relation is equivalent just to the definition (\ref{eq:hatG(x<0,y,m)}), and the second line is justified from the piecewise linearity, connectivity of the function $H_{n,0}(t)$, and the positivity of the function at the ends of each piece, i.e., $H_{n,0}(t)\geq 0$ for $t\in\{-n,-n+2,\cdots,n-2, n\}$, which comes from the property $G(x,y,0)\geq0$ for $x,y\in \mathbb Z_{\geq 0}$ (see Appendix \ref{sec:G(x,y=0,n)} and \ref{sec:G(x,y,n=0)}). Therefore, 
\begin{eqnarray}
\mathbb P\left(0\geq X\;\lor\;-Y\geq 0\big|\mathcal F_0\right)&\leq&H_{X+Y,0}(X-Y)
\end{eqnarray}
where we implicitly use relations (\ref{lemma:upper_ass_1}) and (\ref{lemma:upper_ass_2}).

Next we fix an integer $m_0\in \mathbb Z_{\geq1}$, and we suppose that, in the case of $m=m_0-1\in\mathbb Z_{\geq 0}$, relation (\ref{lemma:upper}) holds (a.s.). We can give 
\begin{eqnarray}
&&\mathbb P\left(X_{m_0}-X_0\geq X\;\lor\; -Y\geq X_{m_0}-X_0\big|\mathcal F_0\right)
\nonumber\\
&=&
\mathbb E\left( \mathbb P\left(X_{m_0}-X_1\geq (X-X_1+X_0)\right.\right.
\nonumber\\&&{}\quad\quad
\left.\left.\;\lor\; -(Y+X_1-X_0)\geq X_{m_0}-X_1\big|\mathcal F_1\right)\big|\mathcal F_0\right)
\nonumber\\
&\leq&\mathbb E\left( H_{X+Y,m_0-1}(X-Y-2X_1+2X_0)\big|\mathcal F_0\right) \quad \makebox{a.s.}
\end{eqnarray}
The first equality comes from the relation $\mathbb P\left(\cdot\big |\mathcal G\right)=\mathbb E\left(\mathbb P\left(\cdot\big|\mathcal G'\right)\big|\mathcal G\right)$ for $\mathcal G\subset \mathcal G'$. To prove the inequality in the third line, we use relation (\ref{lemma:upper}) for $m=m_0-1$ with the substitution $X_n\leftarrow X_{n+1}$, $\mathcal F_n \leftarrow \mathcal F_{n+1}$ for $n\in\{0,1,\cdots \}$, and $X\leftarrow X-X_1+X_0$, $Y\leftarrow Y+X_1-X_0$. Note that, relations (\ref{lemma:upper_ass_1}) and (\ref{lemma:upper_ass_2}) for replaced random variables are given just from the relation $\mathbb E\left(A\big|\mathcal G\right)=A$ $\Rightarrow$ $\mathbb E\left(A\big|\mathcal G'\right)=A$ for $\mathcal G\subset \mathcal G'$. As a result, it is enough to check the relation
\begin{eqnarray}
\mathbb E\left( H_{X+Y,m_0-1}(X-Y-2X_1+2X_0)\big|\mathcal F_0\right)
&\leq  &H_{X+Y,m_0}(X-Y)\quad\makebox{a.s.}
\label{eq:suf_con_0}
\nonumber\\
\end{eqnarray}

A sufficient condition for the above relation is that 
\begin{eqnarray}
\mathbb E\left( H_{n,m-1}(t-2Z)\right)&\leq  &H_{n,m}(t)
\label{eq:suf_con}
\end{eqnarray}
for any non-negative integers $n,m-1\in\mathbb Z_{\geq 0}$, any real number $t\in\mathbb R$, and any random variable $Z$ whose expected value is $0$ and whose absolute value is not more than $1$, i.e., $\mathbb E\left(Z\right)=0$, and $\left|Z\right|\leq 1$. Relation (\ref{eq:suf_con_0}) is derived from relation (\ref{eq:suf_con}) by substituting $\mathbb E\left(\cdot\big|\mathcal F_0\right)$, $X_1-X_0$, $X+Y$, and $X-Y$ into $\mathbb E\left(\cdot\right)$, $Z$, $n$, and $t$, respectively. The substitution of random variables into constants is justified because the variables are fixed numbers under the condition identified by $\mathcal F_0$, i.e., relations (\ref{lemma:upper_ass_1}) and (\ref{lemma:upper_ass_2}) hold. The conditions for $Z$ in this case, i.e., $\mathbb E\left(X_1-X_0\big|\mathcal F_0\right)=0$ and $\left|X_1-X_0\right|\leq1$, are given from assumptions (\ref{eq:main_assumption_2}), (\ref{eq:main_assumption_3}), and (\ref{eq:main_assumption_1}) for $c=1$ (a.s.).

In the following, we show sufficient condition (\ref{eq:suf_con}) by dividing this situation into four cases as follows:

1) When $|t|\geq n$ or $n\in\{0,1\}$ holds,
\begin{eqnarray}
\mathbb E\left( H_{n,m-1}(t-2Z)\right)&\leq  & \mathbb E\left(1\right)
=1=H_{n,m}(t).
\end{eqnarray}
The first equality comes from the property $ H_{n,m}(t')\leq 1$ for $t'\in\mathbb R$ (see Appendix \ref{sec:hatG_bound}). The last relation is just definition (\ref{eq:hatG(x<0,y,m)}) in the case of $|t|\geq n$ or $n=0$, and it is justified in the case of $n=1$ from relation (\ref{eq:hatG(x>0,y>0,m)}) and $G(1,0,m)=G(0,1,m)=1$ for $m\in\mathbb Z_{\geq 1}$ (see Appendix \ref{sec:G(x,y=0,n)}).

2) When $|t|\leq n-2$,
\begin{eqnarray}
&&
\mathbb E\left( H_{n,m-1}(t-2Z)\right)
\nonumber\\
&\leq&
\mathbb E\bigl(\frac{1+Z}2 H_{n,m-1}(t-2)+\frac{1-Z}2 H_{n,m-1}(t+2)\bigr)
\nonumber\\
&=&
 \frac12 H_{n,m-1}(t-2)+\frac12 H_{n,m-1}(t+2)
\nonumber\\
&=&H_{n,m}(t).
\end{eqnarray}
The first inequality comes from the facts that $H_{n,m-1}(t')$ is convex as a function of $t'\in \bigl\{t'\big||t'|\leq n\bigr\}$ (see Appendix \ref{sec:hatG_conv}) and the assumption $\left|Z\right|\leq 1$. The second relation just comes from $\mathbb E\left(Z\right)=0$. The last relation can be derived from the definition of $H$ (see Appendix \ref{sec:hatG_rec}).

3) When $n-2< t < n$ and $n\in \mathbb Z_{\geq 2}$ hold,
\begin{eqnarray}
&&
\mathbb E\left( H_{n,m-1}(t-2Z)\right)
\nonumber\\
&\leq&
\mathbb E\bigl(\frac{1-Z}{1+\frac{n-t}2}+\frac{\frac{n-t}2+Z}{1+\frac{n-t}2}  H_{n,m-1}(t-2)
\bigr)
\nonumber\\
&=&
\frac1{1+\frac{n-t}2}+ \frac{\frac{n-t}2}{1+\frac{n-t}2}H_{n,m-1}(t-2))
\nonumber\\
&\leq&H_{n,m}(t).
\label{eq:tmp_1}
\end{eqnarray}
The first inequality is justified since the relation 
\begin{eqnarray}
H_{n,m-1}(t-2Z)
&\leq&
\frac{1-Z}{1+\frac{n-t}2}+\frac{\frac{n-t}2+Z}{1+\frac{n-t}2}  H_{n,m-1}(t-2)
\end{eqnarray}
holds when $|Z|\leq 1$. This relation itself is guaranteed in the case of $-\frac{n-t}2\leq Z\leq 1$ because of the convexity of $H_{n,m-1}(t')$ for $t-2 \leq t'\leq n$. Even in the case of $-1 \leq Z\leq -\frac{n-t}2$, the left hand side of the relation is equal to $1 = \frac{1-Z}{1+\frac{n-t}2}+\frac{\frac{n-t}2+Z}{1+\frac{n-t}2}$ by definition (\ref{eq:hatG(x<0,y,m)}), and the right hand side of the relation is lower bounded by $1$, which is checked from the non-positivity of $\frac{\frac{n-t}2+Z}{1+\frac{n-t}2} $ and the relation $H_{n,m-1}(t-2)\leq 1$ (see Appendix \ref{sec:hatG_bound}). The second relation of Eq. (\ref{eq:tmp_1}) just comes from $\mathbb E\left(Z\right)=0$. The last relation can be derived from the definition of $H$ (see Appendix \ref{sec:hatG_rec}). 

4) When $-n< t < -n+2$ and $ n\in\mathbb Z_{\geq 2} $ hold, relation (\ref{eq:suf_con}) is directly derived from the third case and the symmetry $H_{n,m'}(t')=H_{n,m'}(-t')$, which comes from the symmetry $G(x,y,m')=G(y,x,m')$, i.e., 
\begin{eqnarray}
\mathbb P\left(H_{n,m-1}(t-2Z)\right)
&=&
\mathbb P\left(H_{n,m-1}(-t-2(-Z))\right)
\nonumber\\
&\leq& H_{n,m}(-t)=H_{n,m}(t).
\end{eqnarray}
The first and the last equalities come from the symmetry, and the second relation can be derived from the inequality in the third case by replacing $t$ and $Z$ with $-t$ and $-Z$. 

\section{Proof of lemma 2}
\label{sec:lem2}
We consider $\{Y_{n+1}^{(x,y)}\}_{x,y,n\in \mathbb Z_{\geq 0}}$ as the random variables, which are defined in lemma 2. From the definition, any probability for each configuration of $\{Y_{n+1}^{(x,y)}\}_{x,y,n\in \mathbb Z_{\geq 0}}$ is uniquely defined. Therefore, we can define 
\begin{eqnarray}
\bar G(x,y,m)&:=&
\mathbb P\bigl(\sum_{n=1}^{m} Y_n^{(x,y)}\geq x\;\lor \;-y\geq \sum_{n=1}^{m} Y_n^{(x,y)}\bigr)
\end{eqnarray}
for any $x,y,m\in\mathbb Z_{\geq 0}$. 

Now, we suppose that $x,y$ are non-negative fixed integers, $X_m$ is a random variable equal to $\sum_{n=1}^{m}Y^{(x,y)}_{n}$, and $\{\mathcal F_n\}_{n\in\mathbb Z_{\geq 0}}$ is a filtration such that $\mathcal F_n:=\mathcal F_n^{(x,y)}$ is a $\sigma$-algebra generated by random variables $Y^{(x,y)}_1, Y^{(x,y)}_2,\cdots, Y^{(x,y)}_n$ for any $n\in\mathbb Z_{\geq 1}$, and $\mathcal F_0=\{\phi,\Omega\}$. From these assumptions, condition (\ref{eq:main_assumption_3}) is obtained. We can also easily check that relations (\ref{eq:main_assumption_1}) and (\ref{eq:main_assumption_2}) for $c=1$ are satisfied since the following two facts can be checked: First, $|Y^{(x,y)}_{m+1}|$ is always not more than $1$. Second, the expected value of $Y^{(x,y)}_{m+1}$ is equal to zero under the condition that only random variables $\{Y^{(x,y)}_{n}\}_{n\in \{1,2,\cdots ,m\}}$ are given. Therefore, the last thing we have to do in the rest of this section is to prove the relation $\bar G(x,y,m)=G(x,y,m)$ for $x,y,m\in \mathbb Z_{\geq 0}$. The structure of the proof is as follows: We give boundary conditions and a recurrence relation which identify $\bar G(x,y,m)$, and we check that $G(x,y,m)$ satisfies the same boundary conditions and the recurrence relation. 

We can evaluate the boundary conditions of the function $\bar G(x,y,m)$ as follows: We can check that both $Y_{n+1}^{(x,0)}$ and $Y_{n+1}^{(0,y)}$ are equal to $0$ with probability $1$ for any $x,y,n\in\mathbb Z_{\geq 0}$. This fact indicates
\begin{eqnarray}
\bar G(x,0,m)&=&\bar G(0,y,m)=\mathbb P\left(0\geq  0\;\lor\;-y\geq 0\right)=\mathbb P\left(0\geq  x\;\lor\;0\geq 0\right)=1
\nonumber\\
\label{eq:opt_bound_1}
\end{eqnarray}
for $x,y,m\in \mathbb Z_{\geq 0}$. The other boundary condition 
\begin{eqnarray}
\bar G(x,y,0)&=&\mathbb P\left(0\geq  x\;\lor\;-y\geq 0\right)=0.
\label{eq:opt_bound_2}
\end{eqnarray}
for $x,y\in \mathbb Z_{\geq1}$ can be derived trivially.

The recurrence relation we use is
\begin{eqnarray}
&&
\bar G(x,y,m) 
\nonumber\\
&=&
\mathbb E\bigl(
\mathbb P\bigl(\sum_{n=2}^m Y_n^{(x,y)}\geq x-Y_1^{(x,y)}\;\lor\;-(y+Y_1^{(x,y)})\geq \sum_{n=2}^m Y_n^{(x,y)}\big|\mathcal F_1^{(x,y)}\bigr)
\bigr)
\nonumber\\
&=&
\mathbb E\bigl(
\bar G(x-Y_1^{(y,z)},y+Y_1^{(y,z)},m-1)
\bigr)
\nonumber\\
&=&
\frac12(\bar G(x-1,y+1,m-1)+\bar G(x+1,y-1,m-1))
\label{eq:opt_rec}
\end{eqnarray}
where $x,y,m \in \mathbb Z_{\geq1}$. The first equality comes from the definition of $\bar G(x,y,m) $ and the property $\mathbb P\left(\cdot\right)=\mathbb E\left(\mathbb P\left(\cdot\big|\mathcal F'\right)\right)$ for any $\sigma$-algebra $\mathcal F'\subset \mathcal F$. In the second equality, we use the relation
\begin{eqnarray}
&&\mathbb P\bigl(Y_1^{(x,y)}=\alpha_1,Y_2^{(x,y)}=\alpha_2,\cdots,Y_{m}^{(x,y)}=\alpha_{m}\bigr)
\nonumber\\
&=&\mathbb P\bigl(Y_1^{(x,y)}=\alpha_1\bigr)\mathbb P\bigl(Y_1^{(x-\alpha_1,y+\alpha_1)}=\alpha_2,\cdots,Y_{m-1}^{(x-\alpha_1,y+\alpha_1)}=\alpha_{m}\bigr)
\end{eqnarray}
which comes from the definition of random variables $Y_n^{(x,y)}$. The above relation guarantees that the arguments of the expectation functions 
 in both sides of the second equation in Eq. (\ref{eq:opt_rec}) are equivalent to each other. The last equality comes from the definition of random variable $Y_1^{(x,y)}$, i.e., when $x$ and $y$ are natural numbers, $Y_1^{(x,y)}$ is equal to $\pm1$ with probability $\frac12$. 

We can easily check that boundary conditions $(\ref{eq:opt_bound_1})$, $(\ref{eq:opt_bound_2})$ and recurrence relation (\ref{eq:opt_rec}) uniquely define the value $\bar G(x,y,m)$ for any $x,y,m\in\mathbb Z_{\geq 0}$. It is straightforward but lengthy to check that $G(x,y,m)$ satisfies the same boundary conditions and recurrence relation. Therefore, we show the proof in an appendix (see Appendix \ref{sec:G(x,y=0,n)}, \ref{sec:G(x,y,n=0)}, and \ref{sec:G_rec}). Combining these results, we have shown the relation $\bar G(x,y,m)=G(x,y,m)$. 

\section{Asymptotics}
\label{sec:comp}
We compare our results with the known bound given from the Azuma-Hoeffding inequality and the bound in the case of a random walk. In this section, we fix parameter $r$ to a positive real number.

We consider the case where $\{X_m\}_{m\in\mathbb Z_{\geq 0}}$ is a martingale that satisfies bounded condition (\ref{eq:main_assumption_1}). From the Azuma-Hoeffding inequality, the following bounds are given:
\begin{eqnarray}
\mathbb P\left(X_m-X_0\geq cr\sqrt{ m}\right)&\leq& \exp(-\frac{r^2}{2}),
\nonumber\\
 \mathbb P\left(|X_m-X_0|\geq cr\sqrt{ m}\right)&\leq&2 \exp(-\frac{r^2}{2}).
\end{eqnarray}
Our tight upper bound gives another asymptotic behavior as follows:
\begin{eqnarray}
&&
 \lim_{m\rightarrow \infty}
\mathbb P\left(X_m-X_0\geq cr\sqrt m\right)
\nonumber\\
&\leq &
\lim_{m\rightarrow \infty} \lim_{y\rightarrow \infty} G(\lfloor r\sqrt m\rfloor,y,m)
\nonumber\\
&=&\lim_{m\rightarrow \infty}
2I_b\bigl(\lfloor \frac{m-\lfloor r\sqrt m\rfloor}2\rfloor,2\lfloor \frac{m-\lfloor r\sqrt m\rfloor}2\rfloor +\lfloor r\sqrt m\rfloor+2\bigr)
\nonumber\\
&=&{\rm erfc}(\frac r{\sqrt2}),
\end{eqnarray}
and 
\begin{eqnarray}
&&
\lim_{m\rightarrow \infty}
\mathbb P\left(|X_m-X_0|\geq cr\sqrt m\right)
\nonumber\\
&\leq &
\lim_{m\rightarrow \infty}
 G(\lfloor r\sqrt m\rfloor,\lfloor r\sqrt m\rfloor,m)
\nonumber\\
&=&
\lim_{m\rightarrow \infty}
\sum_{w\in\mathbb Z_{\geq 0}}
(-1)^w4I_b\bigl(\lfloor \frac{m-\lfloor r\sqrt m\rfloor}2\rfloor -\lfloor r\sqrt m\rfloor w
\nonumber\\&&{}\quad\quad\quad\quad\quad\quad\quad\quad\quad
,2\lfloor \frac{m-\lfloor r\sqrt m\rfloor}2\rfloor +\lfloor r\sqrt m\rfloor+2\bigr)
\nonumber\\
&=&\sum_{w\in\mathbb Z}(-1)^w2\;{\rm erfc}(\frac {(1+2w)r}{\sqrt2})
\end{eqnarray}
where ${\rm erfc}(z)$ is the complementary error function, i.e.,
\begin{eqnarray}
{\rm erfc}(z)&:=&\frac2{\sqrt \pi}\int_{z}^{\infty}e^{-t^2}dt.
\end{eqnarray}

We also discuss the simple bound we give as a corollary. In the case of the one-sided tail probability $\mathbb P\left(X_m-X_0\geq cr\sqrt m\right)$, we can check that the bound gives the same limit value as that for the tight bound. However, in the case of the two-sided probability $\mathbb P\left(X_m-X_0\geq cr\sqrt m\right)$, the bound gives a different limit value as
\begin{eqnarray}
&&
\lim_{m\rightarrow \infty}
\mathbb P\left(\left|X_m-X_0\right|\geq cr\sqrt m\right)
\nonumber\\
&\leq &
\lim_{m\rightarrow \infty}
4I_b\bigl(\lfloor \frac{m-\lfloor r\sqrt m\rfloor}2\rfloor,2\lfloor \frac{m-\lfloor r\sqrt m\rfloor}2\rfloor +\lfloor r\sqrt m\rfloor+2\bigr)
\nonumber\\
&=&2{\rm erfc}(\frac r{\sqrt2}).
\end{eqnarray}

As a reference, we consider the case of a random walk as a special case of a martingale: $\{X_m''\}_{m\in\mathbb Z_{\geq 0}}$ are random variables of a random walk whose one step distance is $c$, i.e.,
\begin{eqnarray}
&&
\mathbb P\left(X_m''= X_{m-1}''+z\big| X_0'', X_1'',\cdots,  X_{m-1}''\right)
\nonumber\\
&=&
\left\{
\begin{array}{ll}
\frac12 & \makebox{in the case of $z\in\{c,-c\}$},
\\
0&\makebox{in other cases}.
\end{array}
\right.
\end{eqnarray}
In this case, we can explicitly derive the tail probabilities as
\begin{eqnarray}
 \lim_{m\rightarrow \infty}
\mathbb P\left(X_m''- X_0''\geq cr\sqrt m\right)&=&
\lim_{m\rightarrow \infty}2^{-m}
\sum_{m=0}^{\lfloor \frac{m- r\sqrt m}2\rfloor }\frac{m!}{n!(m-n)!}
\nonumber\\
&=&\frac12 {\rm erfc}(\frac r{\sqrt2}),
\nonumber\\
 \lim_{m\rightarrow \infty}
\mathbb P\left(\left| X_m''- X_0''\right|\geq cr\sqrt m\right)
&=& 
{\rm erfc}(\frac r{\sqrt2}).
\end{eqnarray}

\begin{figure}
 \begin{center}
  \includegraphics[width=80mm]{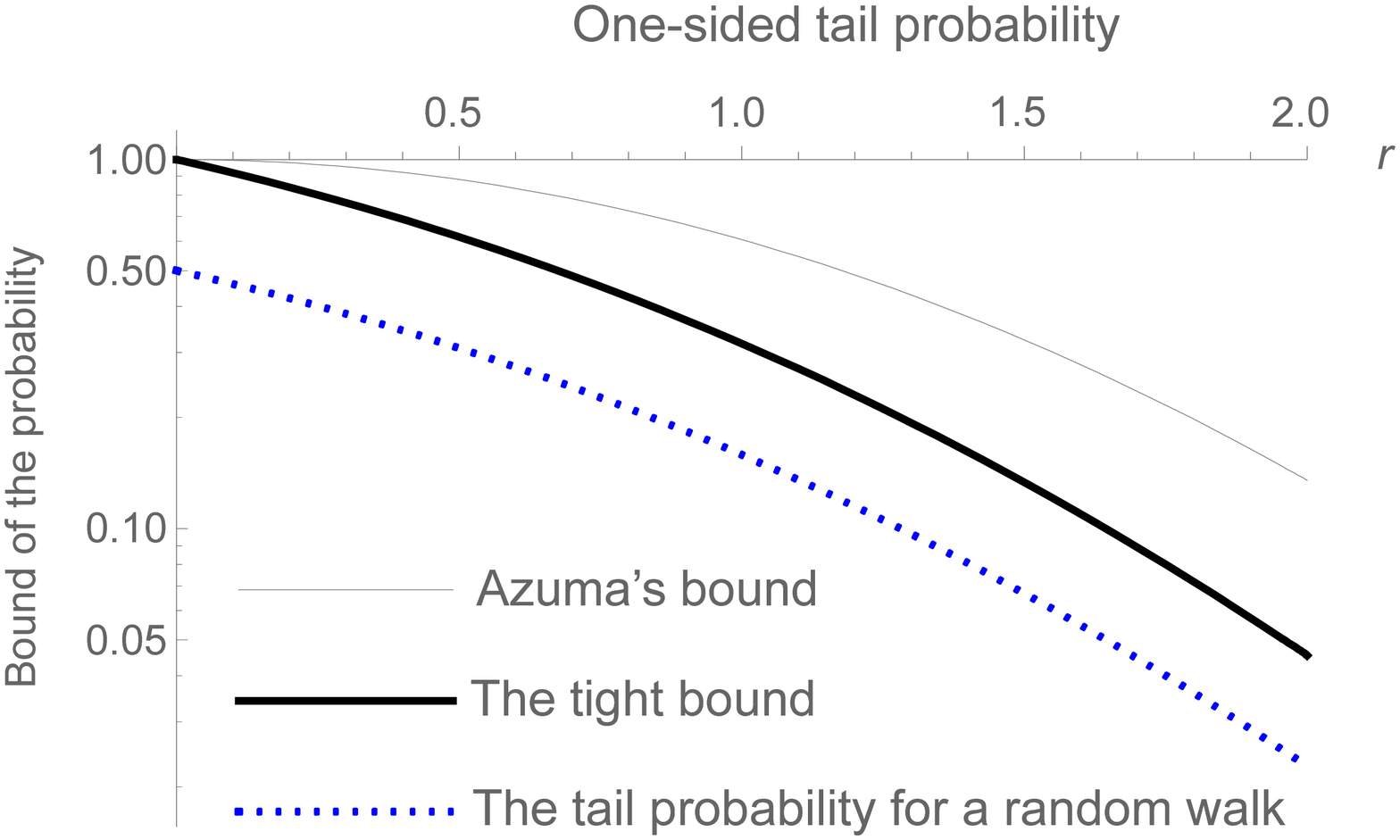}
  \includegraphics[width=80mm]{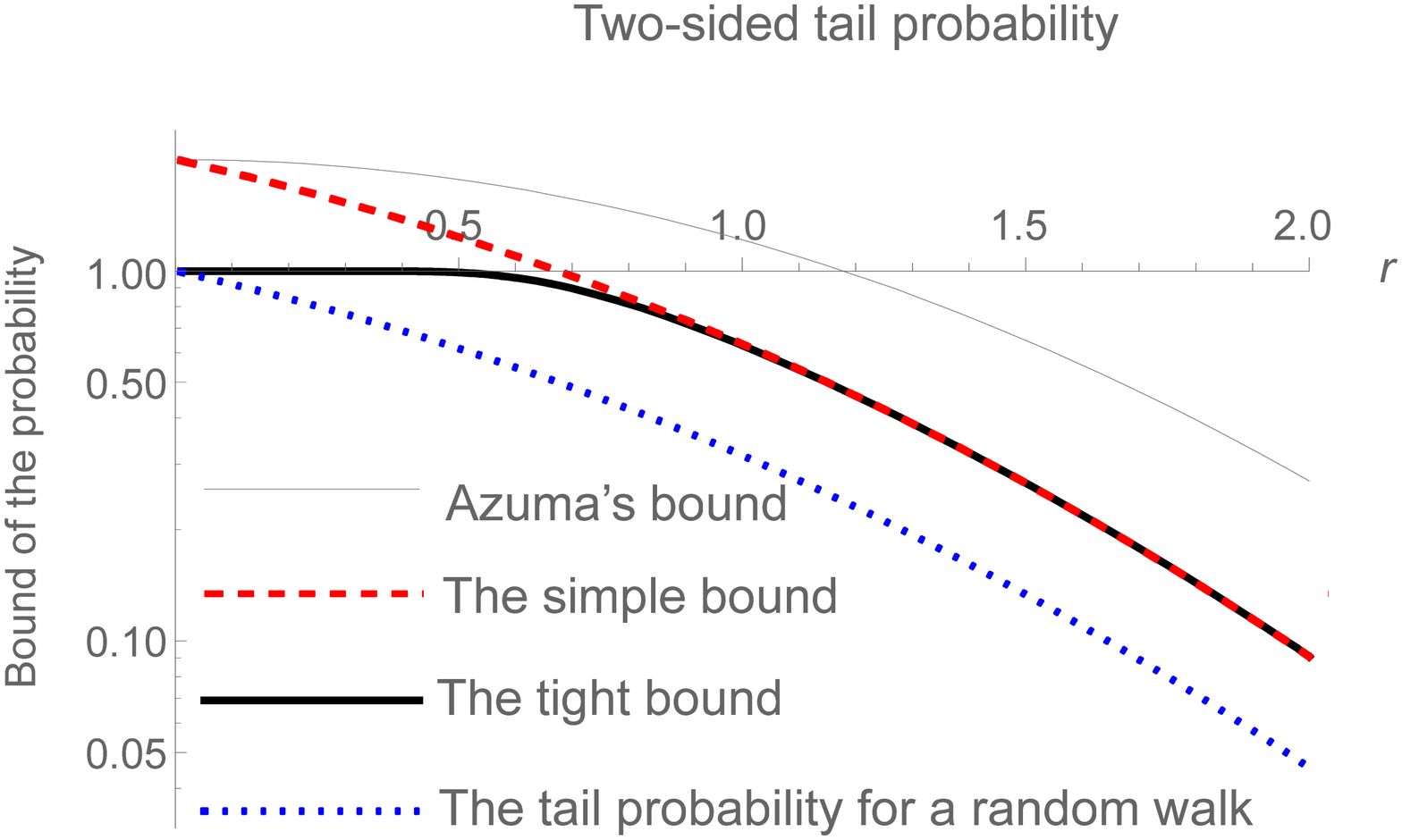}
 \end{center}
\caption{\textbf{ Upper bounds of tail probabilities $\lim_{m\rightarrow \infty} \mathbb P\left(X_m-X_0\geq cr\sqrt m\right)$ and $\lim_{m\rightarrow \infty} \mathbb P\left(\left|X_m-X_0\right|\geq cr\sqrt m\right)$ for a bounded martingale.} The thin gray lines are the bounds derived by the Azuma-Hoeffding inequality. The thick black lines are the tight bounds derived in this paper. The red dashed line indicates the case of the simple bound we give as a corollary. Note that for the one-sided tail probability, the simple bound is equal to the tight bound. The dotted blue lines are included as references. Those are the tail probabilities $\lim_{m\rightarrow \infty} \mathbb P\left(X_m''-X_0''\geq cr\sqrt m\right)$ and $\lim_{m\rightarrow \infty} \mathbb P\left(\left|X_m''-X_0''\right|\geq cr\sqrt m\right)$ in the case that the bounded martingale is a series of random variables generated by a random walk whose one step distance is $c$. Note that $c$ is the parameter which identifies condition (\ref{eq:main_assumption_1}), and the bound doesn't depend on it.\label{fig:bound}}
\end{figure}

All these bounds and the probabilities are shown in Fig.\ref{fig:bound}.

\section{Conclusion}
We gave tight upper bounds of tail probabilities for a bounded martingale. We believe that this result will have benefits in the many fields where the Azuma-Hoeffding inequality is used. Other than such a pragmatic benefit, we hope that the strategy we have used to check the tightness is applicable to derive other concentration inequalities that give tight bounds in other cases.

\label{sec:super_sub}

\appendix

\section{Properties of $G(x,y,m)$}
\subsection{Values of $G(x,y=0,m)$ and $G(x=0,y,m)$}
\label{sec:G(x,y=0,n)}
We show that $G(x,0,m)=G(0,x,m)=1$ for $x,m\in\mathbb Z_{\geq 0}$.

In a special case $x=0$, $G(0,0,m)$ is equal to $1$ by definition.

In other cases, i.e., $x\in\mathbb Z_{\geq 1}$, $G(x,0,m)=G(0,x,m)=1$ can be shown as follows:
\begin{eqnarray}
&&G(x,0,m)=G(0,x,m)
\nonumber\\
&=&
\sum_{w=0}^{\lfloor\frac{m}{2x}\rfloor}
2 I_b\bigl(\lfloor \frac{m}2\rfloor -xw,2\lfloor \frac{m}2\rfloor +2\bigr)
-2 I_b\bigl(\lfloor \frac{m}2\rfloor-x -xw,2\lfloor \frac{m}2\rfloor +2\bigr)
\nonumber\\&&\quad\quad{}
+ 2I_b\bigl(\lfloor \frac{m-x}2\rfloor -xw,2\lfloor \frac{m-x}2\rfloor +x+2\bigr)
\nonumber\\&&\quad\quad{}
-2 I_b\bigl(\lfloor \frac{m-x}2\rfloor -xw,2\lfloor \frac{m-x}2\rfloor +x+2\bigr)
\nonumber\\
&=&2 I_b\bigl(\lfloor \frac{m}2\rfloor,2\lfloor \frac{m}2\rfloor +2\bigr)
=1.
\nonumber\\
\end{eqnarray}
The first relation just comes from the symmetry embedded in definition (\ref{eq:main_func}). The equality in the second line is given just by substituting  definition (\ref{eq:main_func}). In the third equality, we just delete canceled terms, e.g., the third and fourth terms are canceled since they have the same absolute value and the opposite sign, and we also use the fact that $I_b\bigl(\lfloor\frac m2\rfloor -x-x\lfloor \frac m{2x}\rfloor,2\lfloor\frac m2\rfloor+2 \bigr)=0$ for $x\in \mathbb Z_{\geq1}$, which comes from $\lfloor\frac m2\rfloor -x-x\lfloor \frac m{2x}\rfloor<0$. The last equality is derived as follows: For $n\in \mathbb Z_{\geq 0}$, 
\begin{eqnarray}
I_b\bigl(n,2n+2\bigr)
&=&
\sum_{z=0}^n
\frac{1+(-1)^{n-z}}{2^{2n+2}}\frac{(2n+2)!}{z!(2n+2-z)!}
\nonumber\\
&=&
\frac12
\sum_{z=0}^n
\frac{1+(-1)^{n-z}}{2^{2n+2}}\frac{(2n+2)!}{z!(2n+2-z)!}
\nonumber\\&&{}
+
\frac12
\sum_{z=n+2}^{2n+2}
\frac{1+(-1)^{n-z}}{2^{2n+2}}\frac{(2n+2)!}{z!(2n+2-z)!}
\nonumber\\
&=&2^{-2n-3}
\sum_{z=0}^{2n+2}(1+(-1)^{z-n})\frac{(2n+2)!}{z!(2n+2-z)!}
\nonumber\\
&=&2^{-2n-3}
((1+1)^{2n+2}+(-1)^{-n}(1+(-1))^{2n+2})=\frac12.
\end{eqnarray}
The first relation comes from just definition (\ref{eq:qumlant}). In the second equality, we use the fact that the summand in the left hand side as a function of $z$ has a symmetry such that the value does not change for the substitution $z\leftarrow 2n+2-z$. In the third equality, we use the fact that the summand is $0$ when $z=n+1$. The fourth equality is justified simply from a binomial expansion.

\subsection{Value of $G(x,y,m=0)$}
\label{sec:G(x,y,n=0)}
We show that the value $G(x,y,m=0)$ is equal to $0$ for $x,y\in\mathbb Z_{\geq1}$.

The derivation is straightforward as follows:
\begin{eqnarray}
G(x,y,0)
&=&
2 I_b\bigl(\lfloor \frac{-x}2\rfloor,2\lfloor \frac{-x}2\rfloor +x+2 \bigr)
\nonumber\\&&{}
-2 I_b\bigl(\lfloor \frac{-x}2\rfloor-y,2\lfloor \frac{-x}2\rfloor +x+2 \bigr)
\nonumber\\&&{}
+2 I_b\bigl(\lfloor \frac{-y}2\rfloor,2\lfloor \frac{-y}2\rfloor +y+2 \bigr)
\nonumber\\&&{}
-2 I_b\bigl(\lfloor \frac{-y}2\rfloor-x,2\lfloor \frac{-y}2\rfloor +y+2 \bigr)
\nonumber\\
&=&0.
\end{eqnarray}
The first relation comes from just definition (\ref{eq:main_func}). The second equality comes from the fact that $I_b\bigl(-n,m\bigr)=0$ for $m,n\in \mathbb Z_{\geq1}$, which is directly given from the definition of $I_b$.

\subsection{Recurrence relation}
\label{sec:G_rec}
We show the relation 
\begin{eqnarray}
G(x,y,m)&=&\frac12\bigl(G(x-1,y+1,m-1)+G(x+1,y-1,m-1)\bigr)
\label{eq:G_rec}
\end{eqnarray}
for $x,y,m\in \mathbb Z_{\geq1}$.

The derivation is straightforward as follows:
\begin{eqnarray}
&&\frac12(G(x-1,y+1,m-1)+G(x+1,y-1,m-1))
\nonumber\\
&=&
\sum_{w\in\mathbb Z_{\geq0}}
   I_b\bigl(\lfloor \frac{m-x}2\rfloor -(x+y)w,2\lfloor \frac{m-x}2\rfloor +x+1\bigr)
\nonumber\\&&\quad\quad\quad{}
+ I_b\bigl(\lfloor \frac{m-x}2\rfloor-1 -(x+y)w,2\lfloor \frac{m-x}2\rfloor +x+1\bigr)
\nonumber\\&&\quad\quad\quad{}
- I_b\bigl(\lfloor \frac{m-x}2\rfloor-y -(x+y)w,2\lfloor \frac{m-x}2\rfloor +x+1\bigr)
\nonumber\\&&\quad\quad\quad{}
- I_b\bigl(\lfloor \frac{m-x}2\rfloor-y-1 -(x+y)w,2\lfloor \frac{m-x}2\rfloor +x+1\bigr)
\nonumber\\&&\quad\quad\quad{}
+ I_b\bigl(\lfloor \frac{m-y}2\rfloor -(x+y)w,2\lfloor \frac{m-y}2\rfloor +y+1\bigr)
\nonumber\\&&\quad\quad\quad{}
+ I_b\bigl(\lfloor \frac{m-y}2\rfloor-1 -(x+y)w,2\lfloor \frac{m-y}2\rfloor +y+1\bigr)
\nonumber\\&&\quad\quad\quad{}
- I_b\bigl(\lfloor \frac{m-y}2\rfloor-x -(x+y)w,2\lfloor \frac{m-y}2\rfloor +y+1\bigr)
\nonumber\\&&\quad\quad\quad{}
- I_b\bigl(\lfloor \frac{m-y}2\rfloor-x-1 -(x+y)w,2\lfloor \frac{m-y}2\rfloor +y+1\bigr)
\nonumber
\\
&=&
\sum_{w\in\mathbb Z_{\geq0}}
 2I_b\bigl(\lfloor \frac{m-x}2\rfloor -(x+y)w,2\lfloor \frac{m-x}2\rfloor +x+2\bigr)
\nonumber\\&&\quad\quad\quad{}
-2 I_b\bigl(\lfloor \frac{m-x}2\rfloor-y -(x+y)w,2\lfloor \frac{m-x}2\rfloor +x+2\bigr)
\nonumber\\&&\quad\quad\quad{}
+2 I_b\bigl(\lfloor \frac{m-y}2\rfloor -(x+y)w,2\lfloor \frac{m-y}2\rfloor +y+2\bigr)
\nonumber\\&&\quad\quad\quad{}
-2 I_b\bigl(\lfloor \frac{m-y}2\rfloor-x -(x+y)w,2\lfloor \frac{m-y}2\rfloor +y+2\bigr)
\nonumber\\
&=&G(x,y,m).
\end{eqnarray}
The first and the last relation comes from just definition (\ref{eq:main_func}). Note that, though the region of the summation with respect to $w$ is enlarged from $\{0,1,\cdots ,\lfloor \frac{m}{2(x+y)}\rfloor\}$ in definition (\ref{eq:main_func}) into $\mathbb Z_{\geq 0}$ in the above expression, the values are not changed since all terms added are equal to zero. In the second relation, we simplify four adjacent pairs in the left hand side, e.g., the first term and second term, by using the relation $I_b\bigl(n,m\bigr)+I_b\bigl(n-1,m\bigr)=2I_b\bigl(n,m+1\bigr)$ for $m\in \mathbb Z_{\geq1}$ and $n\in\mathbb Z$, which will be proved directly:
\begin{eqnarray}
I_b\bigl(n,m\bigr)+I_b\bigl(n-1,m\bigr)
&=&
\sum_{z=0}^n\frac{1+(-1)^{n-z}}{2^m}\frac{m!}{z!(m-z)!}
\nonumber\\&&{}
+
\sum_{z=0}^{n-1}\frac{1+(-1)^{n-1-z}}{2^m}\frac{m!}{z!(m-z)!}
\nonumber\\
&=&
\sum_{z=0}^n\frac{1+(-1)^{n-z}}{2^m}\frac{m+1-z}{m+1}\frac{(m+1)!}{z!(m+1-z)!}
\nonumber\\
&&{}+\sum_{z=1}^{n}\frac{1+(-1)^{n-z}}{2^m}\frac{z}{m+1}\frac{(m+1)!}{z!(m+1-z)!}
\nonumber\\
&=&2\sum_{z=0}^n\frac{1+(-1)^{n-z}}{2^{m+1}}\frac{(m+1)!}{z!(n+1-z)!}
\nonumber\\
&=&2I_b\bigl(n,m+1\bigr).
\end{eqnarray}
The first and last relation come just from definition (\ref{eq:qumlant}). In the second relation, we just modify the expressions of each term of the summations. And, in the third relation, we use the fact that the summand in the second term in the left-hand side is equal to 0 in the case of $z=0$.

\subsection{Convexity}
\label{sec:G_conv}
We show the convexity of $G(n-t,t,m)$ as a discrete function of $t\in\{0,1,\cdots, n\}$ for fixed $n\in\mathbb Z_{\geq1}$ and $m\in\mathbb Z_{\geq 0}$, i.e.
\begin{eqnarray}
 2G(n-t,t,m)&\leq& G(n-t-1,t+1,m)+G(n-t+1,t-1,m)
\label{eq:G_conv}
\end{eqnarray}
for $n\in\mathbb Z_{\geq 2}$, $m\in\mathbb Z_{\geq 0}$, and $t\in\{1,2,\cdots, n-1\}$.

We prove relation (\ref{eq:G_conv}) by mathematical induction with respect to $m$. 

When $m=0$, relation (\ref{eq:G_conv}) is trivial since we can check it directly by
\begin{eqnarray}
G(n-t,t,0)&=&\left\{
\begin{array}{ll}
1&\makebox{in the case of $t=0$ or $t=n$ }\\
0&\makebox{in the case of $t\in\{1,2,\cdots, n-1\}$}
\end{array}
\right.
\end{eqnarray}
(see Appendix \ref{sec:G(x,y=0,n)} and \ref{sec:G(x,y,n=0)}).

Suppose that the function $G(n-t,t,m)$ is a convex function for fixed $m=m_0-1\in\mathbb Z_{\geq 0}$. We divide the parameter region $n$, $m$ into four regions, and prove the convexity for each region:

1) The first region is $n\in\mathbb Z_{\geq4}$, $t\in\{2,3,\cdots, n-2\}$. In this case, we can prove it without considering the boundary of the function $G$ as follows.
\begin{eqnarray}
&&G(n-t-1,t+1,m_0)+G(n-t+1,t-1,m_0)-2G(n-t,t,m_0)
\nonumber\\
&=&
\frac12\bigl(G(n-t-2,t+2,m_0-1)+G(n-t,t,m_0-1)
\nonumber\\&&{}
-2G(n-t-1,t+1,m_0-1)\bigr)
+
\frac12\bigl(G(n-t,t,m_0-1)
\nonumber\\&&{}
+G(n-t+2,t+2,m_0-1)-2G(n-t+1,t-1,m_0-1)\bigr)
\nonumber\\
&\geq&0.
\nonumber\\
\end{eqnarray}
In the first equality, we use recurrence relation (\ref{eq:G_rec}). In the second relation, the convexity of $G$ at the point $m=m_0-1$ is used.

2) The second region is $n\in\mathbb Z_{\geq 3}$, $t=1$. To prove the convexity, we use an explicit value on the boundary as follows.
\begin{eqnarray}
&&G(n-2,2,m_0)+G(n,0,m_0)-2G(n-1,1,m_0)
\nonumber\\
&=&
\frac12\bigl(G(n-3,3,m_0-1)+G(n-1,1,m_0-1)-2G(n-2,2,m_0-1)\bigr)
\nonumber\\&&{}
+\bigl(G(n,0,m_0)-G(n,0,m_0-1)\bigr)
\nonumber\\
&\geq&0.
\nonumber\\
\end{eqnarray}
In the first equality, we use recurrence relation (\ref{eq:G_rec}) for the first and the second terms in the left-hand side. In the second relation, the convexity of $G$ at  point $m=m_0-1$ and the value on the boundary, i.e., $G(n,0,m_0)=G(n,0,m_0-1)=1$ (see Appendix \ref{sec:G(x,y=0,n)}), are used.

3) The third region is $n\in\mathbb Z_{\geq 3}$, $t=n-1$. We use a symmetry of $G$, i.e., $G(x',y',m)=G(y',x',m)$, and the convexity in the second region.
\begin{eqnarray}
&&G(0,n,m_0)+G(2,n-2,m_0)-2G(1,n-1,m_0)
\nonumber\\
&=&G(n,0,m_0)+G(n-2,2,m_0)-2G(n-1,1,m_0)
\geq 0.
\end{eqnarray}

4) The last region is $n=2$, $t=1$. We can evaluate all the values which appear in Eq. (\ref{eq:G_conv}) explicitly as follows: $G(1,0,m_0)=G(0,1,m_0)=1$ (see Appendix \ref{sec:G(x,y=0,n)}), and $G(1,1,m_0)=\bigl(G(0,2,m_0-1)+G(2,0,m_0-1)\bigr)/2=1$ since $m_0\in\mathbb Z_{\geq 1}$ (see Appendix \ref{sec:G(x,y=0,n)} and \ref{sec:G_rec}). Therefore, we can check that relation (\ref{eq:G_conv}) holds in this region.

\section{Properties of $H_{n,m}(t)$}
\subsection{The convexity of $H_{n,m}(t)$}
\label{sec:hatG_conv}
We confirm the convexity of $H_{n,m}(t)$ as a function of $t\in \bigl\{t\big|| t|\leq n\bigr\}$ for fixed $n,m\in \mathbb Z_{\geq0}$. This is trivially given from the convexity of the discrete function $G(\frac12(n+t),\frac12(n-t),m)$ (see Appendix \ref{sec:G_conv}), and the fact that the discrete function is extended to an affine function $H_{n,m}(t)$ on each interval $n'\leq t\leq n'+2$ for $n'\in\{-n,-n-2,\cdots.n-2\}$.

\subsection{An upper bound of $H_{n,m}(t)$}
\label{sec:hatG_bound}
We show that $H_{n,m}(t)$ is bounded by $1$ for $t\in \mathbb R$ and $n,m\in \mathbb Z_{\geq 0}$.

If $|t|\geq n$, $H_{n,m}(t)$ is equal to $1$ by definition. If $|t|< n$, from the convexity of $H_{n,m}(t))$ proved above, we can evaluate the value as
\begin{eqnarray}
H_{n,m}(t)&\leq&\frac{n+t}{2n}H_{n,m}(n)+\frac {n-t}{2n}H_{n,m}(-n)=1.
\end{eqnarray}
The equality comes from the fact $H_{n,m}(n)=H_{n,m}(-n)=G(n,0,m)=G(0,n,m)=1$ (see Appendix \ref{sec:G(x,y=0,n)}).

\subsection{Weak recurrence relation of $H_{n,m}(t)$}
\label{sec:hatG_rec}
In this subsection, we show the following two relations: 
\begin{eqnarray}
H_{n,m}(t)&=&
\frac12\bigl(H_{n,m-1}(t-2)+H_{n,m-1}(t+2)\bigr)
\label{eq:hatG_rec_1}
\end{eqnarray}
for $|t|\leq n-2$ and $n,m\in \mathbb Z_{\geq1}$, and 
\begin{eqnarray}
H_{n,m}(t)\geq \frac{1}{1+\frac{n-t}2}+\frac{\frac{n-t}2}{1+\frac{n-t}2}H_{n,m-1}(t-2)
\label{eq:hatG_rec_2}
\end{eqnarray}
for $n-2 < t< n$ and $n-1,m\in \mathbb Z_{\geq1}$.

First, we show relation (\ref{eq:hatG_rec_1}):
\begin{eqnarray}
&&H_{n,m}(t))
\nonumber\\
&=&
(1-\frac{n+t}2+z)G(z ,n-z,m)+
(\frac{n+t}2-z)G(z+1,n-z-1,m)
\nonumber\\
&=&
\frac12\bigl(
(1-\frac{n+t-2}2+z_-)G(z_- ,n-z_-,m-1)
\nonumber\\&&{}
+(\frac{n+t-2}2-z_-)G(z_-+1,n-z_--1,m-1)
\nonumber\\&&{}
+(1-\frac{n+t+2}2+z_+)G(z_+ ,n-z_+,m-1)
\nonumber\\&&{}
+(\frac{n+t+2}2-z_+)G(z_++1,n-z_+-1,m-1)
\bigr).
\nonumber\\
&=&
\frac12\bigl(H_{n,m-1}(t-2)+H_{n,m-1}(t+2)\bigr)
\nonumber\\
\end{eqnarray}
where $z:=\lfloor\frac{n+t}2\rfloor$ and $z_\pm:=\lfloor\frac{n+t\pm2}2\rfloor=z\pm1$. In the first and the last relations, we use definition (\ref{eq:hatG(x>0,y>0,m)}) of $H_{n,m}(t))$. In the second relation, we use the recurrence relation (\ref{eq:G_rec}) of $G(x,y,t)$.

Next, we show relation (\ref{eq:hatG_rec_2}):
\begin{eqnarray}
&&
H_{n,m}(t)
-\frac{1}{1+\frac{n-t}2}-\frac{\frac{n-t}2}{1+\frac{n-t}2}H_{n,m-1}(t-2)
\nonumber\\
&=&
\frac{n-t}2G(n-1 ,1,m)
-\frac{(\frac{n-t}2)^2}{1+\frac{n-t}2}
\nonumber\\&&{}
-\frac{(\frac{n-t}2)^2}{1+\frac{n-t}2}G(n-2 ,2,m-1)
-\frac{\frac{n-t}2(1-\frac{n-t}2)}{1+\frac{n-t}2}
G(n-1,1,m-1)
\nonumber\\
&=&
\frac{\frac{n-t}2(1-\frac{n-t}2)}{2(1+\frac{n-t}2)}\bigl(
G(n,0,m-1)+G(n-2 ,2,m-1)
\nonumber\\&&{}\quad\quad\quad\quad\quad\quad
-2G(n-1,1,m-1)\bigr)
\nonumber\\
&\geq&0.
\nonumber\\
\end{eqnarray}
In the first relation, we use definition (\ref{eq:hatG(x>0,y>0,m)}) of $H_{n,m}(t)$, and simplify the expressions by using simple relations like $G(n,0,m)=1$ and $\lfloor \frac{n+t}2\rfloor=n-1$ in this situation. The second relation can be obtained by applying the recurrence relation (\ref{eq:G_rec}) for the first term in the left hand side. The last inequality is justified from the positivity of the coefficient and the convexity (\ref{eq:G_conv}) of $G(n-t,t,m-1)$.

\end{document}